\documentclass[12pt]{article}
\usepackage{amssymb,amscd}

\usepackage{amssymb}
\usepackage{graphicx}
\usepackage{amsfonts}
\usepackage{latexsym}
\usepackage{amsthm}
\usepackage{amsmath}
\usepackage{stmaryrd}

\usepackage{amssymb, amscd}

\usepackage{url}
\usepackage{mathtools}
\usepackage{color}

\newcommand{\zz}[1]{}

\newtheorem{conj}{Conjecture}[section]
\newtheorem{theo}[conj]{Theorem}
\newtheorem{rem}[conj]{Remark}

\newtheorem{defin}[conj]{Definition}
\newtheorem{prop}[conj]{Proposition}
\newtheorem{cor}[conj]{Corollary}
\newtheorem{lema}[conj]{Lemma}

\newtheorem{example}[conj]{Example}

\begin{document}
\date{September 15, 2018}
\title{\Large Topology  of unavoidable complexes}

\author{
{Du\v{s}ko Joji\'{c}} \\ {\small Faculty of Science}\\[-2mm] {\small University of Banja Luka}
\and Wac{\l}aw Marzantowicz \\ {\small Faculty of Mathematics and Computer Science}\\[-2mm] {\small Adam Mickiewicz University , Pozna\'{n}}
\and Sini\v{s}a T. Vre\'{c}ica\\ {\small Faculty of Mathematics}\\[-2mm] {\small University of Belgrade}
\and Rade T. \v{Z}ivaljevi\'{c}\\ {\small Mathematical Institute}\\[-2mm] {\small SASA,    Belgrade}\\[-2mm]}


\maketitle 

\begin{abstract}\noindent
The partition number $\pi(K)$ of a simplicial complex $K\subseteq
2^{[m]}$ is the minimum integer $\nu$ such that for each partition
$A_1\uplus\ldots\uplus A_\nu = [m]$ of $[m]$ at least one of the
sets $A_i$ is in $K$. A complex $K$ is {\em $r$-unavoidable} if
$\pi(K)\leq r$. We say that a complex $K$ is {\em almost
$r$-non-embeddable} in $\mathbb{R}^d$ if for each continuous map
$f: \vert K\vert \rightarrow \mathbb{R}^d$ there exist $r$ vertex
disjoint faces $\sigma_1,\ldots, \sigma_r$ of $\vert K\vert$ such
that $f(\sigma_1)\cap\ldots\cap f(\sigma_r)\neq\emptyset$.
Motivated by the problems of Tverberg-Van Kampen-Flores type we
prove several results (Theorems~\ref{thm:izvor}, \ref{thm:oblak},
\ref{thm:Dule-Sinisa-Rade}) which link together the combinatorics
and topology of these two classes of complexes. One of our central
observations (Theorem~\ref{thm:Dule-Sinisa-Rade}), summarizing and
extending results of G.~Schild, B.~Gr\"{u}nbaum and many others,
is that interesting examples of (almost) $r$-non-embeddable
complexes can be found among the joins $K = K_1\ast\ldots\ast K_s$
of $r$-unavoidable complexes.
\end{abstract}

\section{Introduction}

`Unavoidable complexes' were originally introduced as `Tverberg unavoidable subcomplexes'  by Blagojevi\'{c}, Frick, and Ziegler in \cite[Section~4]{bfz}.

\medskip
A systematic study of unavoidable complexes, as ``combinatorial objects that may have some independent interest and which may deserve to be studied in their own right'', was initiated in \cite{bestiary}.

\medskip
In this paper we pave the way for the study of the topology of unavoidable complexes, emphasizing the relationship between the partition invariant $\pi(K)$ of a simplicial complex $K$ on one side, and the equivariant index ${\rm Ind}_G(K^{\ast r}_\Delta)$ and the $G$-genus $\gamma_G(K^{\ast r}_\Delta)$ of the associated deleted join (deleted product), on the other.

  \subsection{The `constraint method' of \cite{bfz} and \cite{g10}}

The method of `Tverberg unavoidable complexes' or the `constraint
method', as introduced and developed in \cite{bfz} by Blagojevi\'
c, Frick, and Ziegler, and earlier (in a much less explicit form) by
Gromov \cite{g10}, has proven to be a powerful and
versatile method for generating statements of Tverberg type.

One of the key contributions of \cite{bfz} was the introduction of `Tverberg unavoidable complexes' (see our Definition~\ref{def:BFZ-unavoidable} in
Section~\ref{sec:r-unavoidable}). This concept and its far
reaching and beautiful applications were our main motivation for
isolating the `partition number' $\pi(K)$ in \cite{bestiary}(see also the unpublished preprint \cite{jvz-3}) and the associated class of $r$-unavoidable complexes, as interesting combinatorial objects that deserve to be studied in their own right.

\subsection{$2$-unavoidable complexes}

{\em Unavoidable} or more precisely {\em $2$-unavoidable
simplicial complexes} can be directly linked to some classes of
simplicial complexes which (often independently) emerged and
attracted attention of researchers in {\em game theory,
combinatorial topology, social choice theory, reliability theory},
{\em geometry of moduli spaces of polygonal linkages}, and other
areas. In topology \cite{M} they appear as the Alexander self-dual
complexes, and provide key examples of $n$-dimensional complexes
non-embeddable in $\mathbb{R}^{2n}$. In social choice theory
(reliability theory) \cite{vNM44, r90, ps07}  they arise as the
complexes describing the winning (losing) coalitions in game
theory (`simple games' of Von~Neumann and Morgenstern
\cite{vNM44}). In geometry of configuration spaces \cite{ga-pa}
they arise as the complexes of `short sets', characterizing the
configuration spaces of polygonal linkages in Euclidean spaces.

Our original motivation was the quest for more `exotic' examples of
theorems of Tverberg-Van Kampen-Flores type. However it may be
expected that the invariant $\pi(K)$ and the associated
$r$-unavoidable complexes are as important and interesting for
other fields as their $2$-unavoidable counterparts.

\subsection{This paper}

Motivated by the problems of Tverberg-Van Kampen-Flores type we
prove several results (Theorems~\ref{thm:izvor}, \ref{thm:oblak},
\ref{thm:Dule-Sinisa-Rade}) which link the combinatorics of
$r$-unavoidable complexes to the topology of almost
$r$-non-embeddable complexes.

\medskip
The main idea is to compare the combinatorial complexity of a
simplicial complex $K$ (expressed in terms of the invariant
$r=\pi(K)$) to the (equivariant) topological complexity of the
associated deleted join $K^{\ast r}_\Delta$  (or deleted product
$K^r_\Delta$) evaluated by an associated equivariant index ${\rm
Ind}_G$ (or alternatively by the $G$-genus $\gamma_G$) . This comparison is illustrated by the inequalities proved in Theorems~\ref{thm:izvor} and \ref{thm:oblak}. For example, a corollary of Theorem~\ref{thm:izvor}
(Corollary~\ref{cor:pi}) says that there is an inequality,
\begin{equation}\label{eqn:izvor-3}
{\rm Ind}_G (K^{\ast r}_\Delta) \geq m - \pi(K)
\end{equation}
where $r = \pi(K)=p^k$ is a prime power, $G=(\mathbb{Z}_p)^k$,
$K\subseteq 2^{[m]}$ and ${\rm Ind}_G$ is the equivariant index function described in Section~\ref{sec:eq-index-glavna}.

\medskip
Our central result is Theorem~\ref{thm:Dule-Sinisa-Rade}. By
summarizing and extending several results of  Gr\"{u}nbaum,
Sarkaria, Schild, Blagojevi\'{c}, Frick, Ziegler, as well as our
own work, we demonstrate that many interesting examples of
almost $r$-non-embeddable complexes can be found among the joins
$K = K_1\ast\ldots\ast K_s$ of $r$-unavoidable complexes.

\medskip
\begin{rem}{\rm
It may be instructive and interesting to compare the inequality (\ref{eqn:izvor-3}) with similar inequalities from \cite{M} where an appropriate topological index provides an upper or lower bound for an interesting combinatorial invariant. For example the ``Sarkaria's coloring/embedding theorem'' is in \cite[Theorem~5.8.2]{M} formulated as the inequality
\[
   {\rm Ind}_{\mathbb{Z}_2}(K^{\ast 2}_\Delta)  \geq n - \chi({\rm KG}(\mathcal{F})) -1
\]
where $\chi({\rm KG}(\mathcal{F}))$ is the chromatic number of the Kneser graph ${\rm KG}(\mathcal{F})$ (\cite[Section~3]{M}) where $\mathcal{F}\subset 2^{[n]}$ is the collection of all minimal non-faces of a simplicial complex $K\subseteq 2^{[n]}$.
}
\end{rem}

\section{The invariant $\pi(K)$ and $r$-unavoidable complexes}
\label{sec:r-unavoidable}

The following definition  (Definition~4.1 in
\cite[Section~4]{bfz}) is central for the `constraint method' or
the method of `Tverberg unavoidable complexes'.

\begin{defin}\label{def:BFZ-unavoidable} {\em (Tverberg unavoidable subcomplexes)}.
Let $r\geq 2, d \geq 1, N \geq r-1$ be integers and $f : \Delta_N
\rightarrow \mathbb{R}^d$ be a continuous map with at least one
Tverberg $r$-partition. Then, a subcomplex $\Sigma \subseteq
\Delta_N$ is {\em Tverberg unavoidable} if for every Tverberg
partition $\{\sigma_1, \ldots, \sigma_r\}$ for $f$, there is at
least one face $\sigma_j$ that lies in $\Sigma$.
\end{defin}

\noindent
Recall that a family $\{\sigma_1, \ldots, \sigma_r\}$ of pair-wise vertex disjoint faces of $\Delta_N$ is a ``Tverberg
partition  for $f$'' if $f(\sigma_1)\cap\dots\cap f(\sigma_r)\neq \emptyset$.

\medskip
As remarked in \cite[Section~4]{bfz}, the property of being
`Tverberg unavoidable' (as introduced in
Definition~\ref{def:BFZ-unavoidable}) depends both on the
parameters $r, d$ and $N$, and on the chosen map $f$. However, the
authors of \cite{bfz} emphasized that their main interest in that
paper were the subcomplexes that are large enough to be
unavoidable for any continuous map $f$.

\medskip
In the following closely related definition (see \cite[Definition~2.4.]{bestiary}) we
 avoid any reference to the continuous map $f :
\Delta_N\rightarrow \mathbb{R}^d$ and in particular put less emphasis on
the parameters $d$ and $N$. As a consequence our definition is conceptually simpler and our class of $r$-unavoidable complexes is
larger than the original class of Tverberg unavoidable complexes.

\medskip We say that a family $\{A_i\}_{i=1}^r\subset 2^{[m]}$ is a {\em subpartition of $[m]$} (and write $A_1\uplus\ldots\uplus A_r \subseteq [m]$) if the sets $A_i$ are pair-wise disjoint.

\begin{defin}{\em ($r$-unavoidable complexes)} \label{def:neizbezno}
Let $r\geq 2$ be an integer. Suppose that $K$ is a simplicial
complex with vertices in $[m]$, meaning that $Vert(K)\subseteq [m]$.
The complex $K$ is called $r$-unavoidable on $[m]$ \mbox{\rm (or
simply $r$-unavoidable)} if,
\begin{equation}\label{eqn:uslov-neizbeznosti}
\forall A_1,\ldots, A_r \in 2^{[m]}\setminus\{\emptyset\}, \quad
A_1\uplus\ldots\uplus A_r \subseteq [m] \quad \Rightarrow \quad
(\exists i) \, A_i\in K.
\end{equation}
\end{defin}

\medskip
If $K\subseteq L\subseteq 2^{[m]}$ and the complex $K$ is
$r$-unavoidable then $L$ is $r$-unavoidable as well. This is the
reason why it may be sometimes useful to focus on `minimal
$r$-unavoidable complexes'.

\begin{defin}\label{def:minimal}
A complex $K\subseteq 2^{[m]}$ is {\em minimally $r$-unavoidable} if
it is $r$-unavoidable and if $L\varsubsetneq K$ is a proper
subcomplex of $K$ then $L$ is not $r$-unavoidable.
\end{defin}

The property of being $r$-unavoidable is an intrinsic,
combinatorial property of the simplicial complex $K$. Here is a
natural generalization.

\begin{defin}{\em ($(r,s)$-unavoidable complexes)} \label{def:r-and-s}
Choose integers $r> s\geq  1$. Suppose that $Vert(K)\subseteq [m]$.
The complex $K$ is called {\em $(r,s)$-unavoidable} if,
\begin{equation}\label{eqn:uslov-r-and-s}
\forall A_1,\ldots, A_r \in 2^{[m]}\setminus\{\emptyset\}, \quad
A_1\uplus\ldots\uplus A_r \subseteq [m] \quad \Rightarrow \quad
\vert K\cap \{A_i\}_{i=1}^r\vert \geq s.
\end{equation}
\end{defin}
In other words the condition (\ref{eqn:uslov-r-and-s}) says that
for each partition $\uplus_{i=1}^r~A_i=[m]$ of $[m]$ into $r$
non-empty sets, at least $s$ of the sets $A_i$ belong to $K$.

\medskip Perhaps the most elegant way to introduce the
$r$-unavoidable complexes is via the `partition invariant'
$\pi(K)$. Note that by assuming $K\subseteq 2^{[m]}$ we mean that
${\rm Vert}(K)\subseteq [m]$ and that this inclusion may be strict
in general.

\begin{defin}\label{def:pi-invariant}
The partition number $\pi(K)$ of a simplicial complex $K\subseteq
2^{[m]}$ is the minimum integer $\nu$ such that for each partition
$A_1\uplus\ldots\uplus A_\nu = [m]$ of $[m]$ at least one of the
sets $A_i$ is in $K$. In other words $K\subseteq
2^{[m]}$ is $r$-unavoidable if and only if $\pi(K)\leq r$.
\end{defin}

\section{Equivariant index and $r$-unavoidable complexes}
\label{sec:eq-index-glavna}

In this section we establish  a connection between $r$-unavoidable
complexes and the equivariant index theory.

\medskip
For the reader's convenience we include the definition and outline
the main properties of the numerical index
function ${\rm Ind}_G$ for the elementary abelian $p$-group
$G = (\mathbb{Z}_p)^k$.
This index function agrees with the index function
${\rm Ind}_G$ described in \cite{M} in the case $G = \mathbb{Z}_p$ and
has the merit to assign `correct values' to some important spaces
(spheres) with fixed point free $(\mathbb{Z}_p)^k$-actions.

\medskip
For a non-specialist interested mainly in applications, the equivariant index theory can
be understood as a kind of complexity theory for $G$-complexes. It
provides a source of `Borsuk-Ulam' type results needed for the
application of the {\em Configuration space/test map scheme}
\cite{Ziv-96-98, M, Z04,Z17} and in this sense it can be used as a `black
box' for immediate applications in discrete geometry and
combinatorics.

\subsection{The equivariant index function
${\rm Ind}_G$}\label{sec:extended}

Our index function is a close relative of the $G$-genus, described and developed in
\cite[Def. 2.8]{Bar}. For our purposes it is sufficient to use the variant of $G$-genus
defined in proposition  \cite[Prop. 2.9]{Bar} where the defining
family $\mathcal{A}$ of $G$-spaces is chosen to be the family of  orbits
$(G/H)$, $H\neq G$ of a given finite group $G$.

\begin{defin}\label{genus}
For a given $G$-space $X$  the $G$-genus $\gamma_G(X)$ of $X$ is defined as the smallest number $k$  such that here exists a $G$
equivariant map $$\phi: X \to G/H_1 \,*\, \cdots\,*\, G/H_k$$ where
$H_i\subsetneq  G $ is a  proper subgroup of $G$.
\end{defin}

\begin{prop}\label{prop:genus-prop}
The $G$-genus has the following properties:

\begin{itemize}\label{properties of genus}
\item[(i)] { $\gamma_G(X)\in \mathbb{N}_0 \cup \infty$, $\gamma(G/H)=1$ for
every $H \varsubsetneq G$, $\gamma_G(X) =0 $ if and only if
$X=\emptyset$, and if $X^G \neq \emptyset$ then
$\gamma_G(X)=\infty$.}

\item[(ii)] {Monotonicity: if $f: X\to Y$ is a $G$-map, then
$\gamma_G(X)\leq \gamma_G(Y)\,.$}

\item[(iii)] {Subadditivity: If $X_1, \; X_2$ are two invariant
subsets of a normal space $X$ such that their interiors cover $X$
then $ \gamma_G(X)\;\leq \; \gamma_G(X_1) + \gamma_G(X_1)\,.$
\newline (Moreover the $G$-genus is maximal among the functions satisfying (i)--(iii).) }

\item[(iv)] {Continuity: Every closed
invariant subset $X$  of a metrizable $G$-space $Z$ (with an
invariant metric) has  an open invariant neighborhood $U\supseteq X$
such that $ \gamma_G(U)= \gamma_G(X)\,. $}

\item[(v)] {Finiteness: If $X$ is a compact $G$-space such that
$X^G=\emptyset$ then
\begin{itemize}
\item[a)] {$\gamma_G(X) < \infty$.}
\item[b)]{If $X= X_{(H)} := \{x\in X \mid G_x  \sim H\}$ is a $G$-space with one orbit type, e.g. if the action is free, then $\gamma_G(X)\leq \dim (X/G) +1 $, or equivalently $\gamma_G(X)\leq \dim (X) +1
$, since $G$ is finite.}
\end{itemize}}
\end{itemize}
\end{prop}

\begin{proof} The proofs of all facts {\em (i)--(v)} are elementary and direct. The details can be found in \cite[Propositions 2.15, 2.16]{Bar}, however  the reader can easily
check them as an exercise. \end{proof}

\begin{example}
{\rm  The following ``join property'' is a formal consequence of {\em (ii)} and {\em (iii)} but it can be also deduced directly from the definition of the $G$-genus.}
 \begin{itemize}
\item[(vi)] {Join Property: If $X\ast Y$ is a join of $G$-spaces then $\;\; \gamma_G(X*Y)\,\leq\, \gamma_G(X) + \gamma_G(Y)$.}
\end{itemize}
\end{example}

\begin{rem}{\rm
 The idea of the $G$-genus goes back
to M. Krasnoselski and C. T. Yang who defined and applied the
$\mathbb{Z}_2$-genus in the early 1950s.  The reader is referred to \cite{Bar} for the
references and related information. }
\end{rem}

\medskip
The following proposition can be directly linked to the so called {\em Sarkaria's inequality}, see Proposition~\ref{prop:ind-sharp}~(6) and \cite{Ziv-96-98,M}, which has found numerous applications in combinatorics and discrete geometry.

\begin{prop}($G$-genus inequality for posets)\label{Sarkaria inequality}
 Suppose that $L_0$ is a finite $G$-simplicial complex and let
$L\subseteq L_0$ be a $G$-invariant subcomplex. Then,
\begin{equation}\label{eqn:Sark-inequality1}
  \gamma_G(L_0) \,\leq \gamma_G(L) +
\gamma_G(\Delta(L_0\setminus L))
\end{equation}
where $\Delta(L_0\setminus L)$ is the order complex of the poset
$(L_0\setminus L, \subseteq)$.
\end{prop}
\begin{proof}
It is easy to observe that there exists a $G$-equivariant map
\[
   \Delta(L_0) \stackrel{G}{\longrightarrow} \Delta(L) \ast
   \Delta(L_0\setminus L).
\]
Then both the relation (\ref{eqn:Sark-inequality1}) and Sarkaria's inequality (Proposition~\ref{prop:Sarkaria}~(\ref{eqn:Sark-inequality})) are deduced from the
 Join Property {\em (vi)} of the $G$-genus $\gamma_G$.
\end{proof}

All the properties {\em (i)--(vi)},  listed so far,  hold unconditionally for any finite group $G$. Now we turn our attention to other important, albeit less general properties of the $G$-genus, where it is essential to work with smaller classes of finite groups.

\begin{prop}\label{genus of sphere}
Let $G$ be  a finite abelian group.   Let $V$ be a complex,
orthogonal $G$-representation of real dimension $\dim_\mathbb{R} V = 2\,
\dim_\mathbb{C}  V =
  n $ such that $V^G=\{0\}$. Then
$\gamma_G(S(V))\,\leq \, n$.
 \end{prop}

\begin{proof}  Let $V={\oplus_\alpha}\, k_\alpha V_\alpha$ be the decomposition of $V$ into irreducible factors.  Since $S(W_1\oplus W_2)\cong  S(W_1)* S(W_2)$  for each two $G$-representations $W_1$ and $W_2$  it is sufficient, in light of the property {\em (vi)}, to prove that $\gamma_G(S(V_\alpha))=2$ for each irreducible (nontrivial) complex representation $V_\alpha$.

To prove this note that for a homomorphism $\rho: G \to U(1)= S^1$ the
image ${\rm Image}(\rho(G))\subseteq  S^1$ is a cyclic subgroup $K \subseteq
S^1$, $K \simeq  G/H$, $H=\ker \rho$,   which can be identified with
the roots of unity of order $m$ dividing $ \vert G\vert $. Moreover
the action of $H$ on $S(V)$ is free.

From here one can construct directly a $G$-map from $S(V)$ into $G/H *
G/H$, or alternatively use the inequality $\gamma_G(S(V)) \leq {\rm
cat}_G(S(V))={\rm cat}_K(S(V)) = {\rm cat}(S(V)/K)= {\rm
cat}(S(V)/G)=2$ (cf. \cite{Bar}).
\end{proof}

\begin{rem}\label{Remark on genus of spheres for p-tori}{\rm
There are two important special cases of Proposition~\ref{genus of sphere} where it suffices to assume that $V$ is a {\em real} $G$-representation.  }

\begin{itemize}
\item If $G= \mathbb{Z}_2^k$, $k\geq 1$ is a $2$-torus then every
irreducible real representation of $G$ is one-dimensional, given by
a homomorphism $\rho: G \to \mathbb{Z}_2 =\{-1, \,1\}\subseteq
\mathbb{R}$. Consequently, by the  argument used above we have
$\gamma_G(S(V)) \leq n$ for every real, $n$-dimensional representation $V$ of $G$ such that $V^{G} = \{0\}$.

\item If $G=\mathbb{Z}_p^k$, $k\geq 1$, $p$ odd prime, then every
representation $V$  of $G$ such that $V^G=\{0\}$ possess a complex
structure. Consequently  Proposition~\ref{genus of sphere} still holds.
\end{itemize}
\end{rem}

For the remaining properties of $G$-genus, especially for the
estimates of $\gamma_G$ from above, we need an even more restrictive condition on $G$. From here on we work with the class of  $p$-tori (elementary abelian $p$-groups) $(\mathbb{Z}_p)^k$, where  $p$ is a prime.  This restriction is essential since these properties are deduced from a version of the Borsuk-Ulam theorem which is known not to  hold beyond this class of groups.

\begin{theo}\label{Borsuk-Ulam}
Assume that $G= \mathbb{Z}_p^k $, $k\geq 1$,  is a $p-$torus.
\begin{itemize}
\item[(a)]{ If $V$ is an orthogonal representation of $G$ of real dimension $n$ with
$V^G=\{0\}$  then  $\gamma_G(S(V))= n\, .$}
\item[(b)]{ If $X$ is a fixed point free $G$-space (i.e.
$X^G=\emptyset$) such that $H^i(X; \mathbb{Z}_p) = 0$ for $1\leq i
\leq n-1$ (for example if $X$ is $n-1$-connected) then $ \gamma_G(X) \geq
n+1 \,.$}
\end{itemize}
\end{theo}

\begin{proof}  The proofs of these results use the Borel localization theorem, and (part (b)) a spectral sequence argument. Standard references are \cite[Thm. 2.9]{Marzantowicz1} and \cite[Prop 1.3]{CP1} (for part (a)) and \cite[Prop 6.3]{CP1} (for part  (b)).
\end{proof}

\begin{rem}{\rm
  Note that part (a) of Theorem~\ref{Borsuk-Ulam} can be also deduced from
part (b) and Proposition~\ref{genus of sphere}.  }
\end{rem}

In combinatorial applications, see for example the references in \cite{Z17, M},    the use of the {\em equivariant index}  ${\rm Ind}_G$ (rather than the $G$-genus) is more customary. Since the  $G$-index is related to the $G$-genus in the same
way as the dimension of sphere is related to the dimension of its
Euclidean space, there is no a real difference between them and using one or the other is largely a matter of taste.

\begin{defin}\label{index}
Let $G$ be a finite group and  $\mathcal{C}_G$ be the category of
 $G$-spaces with $G$-equivariant maps
as morphisms.

The associated index function is  defined on $\mathcal{C}_G$ by the
formula,
\begin{equation}\label{eqn:index}
{\rm Ind}_G(X) \;:= \;\gamma_G(X) -1
\end{equation}
\end{defin}

The following proposition is essentially a translation of known facts about the $G$-genus into the language of index function. Note that the property (6) says that this index function agrees with the function described in \cite{M} in the case when $G = \mathbb{Z}_p$.

\begin{prop}\label{prop:ind-sharp}
The equivariant index has the following properties:
\begin{enumerate}
 \item[{\rm (1)}] \quad ${\rm Ind}_G(X) \in
 \mathbb{N}_0\cup\{+\infty\}$.
 \item[{\rm (2)}] Finiteness:  \quad ${\rm Ind}_G(X)<+\infty$ if $X$ is
 compact and $X^G =\emptyset$, \\ ${\rm Ind}_G(X) \leq  \dim (X/G)= \dim (X)$ if $X=X_{(H)}$, e.g. if the action is free.
 \item[{\rm (3)}] Monotonicity: \quad If $ X
\stackrel{G}{\longrightarrow} Y$ then ${\rm Ind}_G(X)\leq
{\rm Ind}_G(Y)$.
 \item[{\rm (4)}]{Join property and subadditivity: \quad ${\rm Ind}_G(X\ast Y)\leq {\rm Ind}_G(X) + {\rm Ind}_G(Y)
 +1$, \quad  \\ and \quad $ {\rm Ind}_G(X) \leq {\rm
 Ind}_G(X_1) +{\rm Ind}_G(X_2) -1$  \\ if  $X_1, \; X_2$ are two invariant
subsets of a normal space $X$ such that their interiors cover $X$.
 }
 \item[{\rm (5)}] \quad If $X$ is  $(n-1)$-connected then ${\rm Ind}_G(X)\geq
 n$.

\item[{\rm (6)}] \quad If $G= \mathbb{Z}_p$,   $p$ is a prime, then ${\rm
Ind}_G(X) $ is the same as the equivariant index described in \cite{M}.
\end{enumerate}
\end{prop}

\medskip
The most important for applications is the following property of the index function, known in combinatorial circles as the {\em Sarkaria's inequality}  (introduced in \cite{Ziv-96-98}, see also  \cite[page 124]{M}).

\begin{prop}{\em (Sarkaria's inequality)}\label{prop:Sarkaria}
Let $G = (\mathbb{Z}_p)^k$ be a $p$-torus.
Suppose that $L_0$ is a finite $G$-simplicial complex and
$L\subseteq L_0$ its $G$-invariant subcomplex. Then,
\begin{equation}\label{eqn:Sark-inequality}
{\rm Ind}_G(L) \geq {\rm Ind}_G(L_0) - {\rm
Ind}_G(\Delta(L_0\setminus L)) -1,
\end{equation}
where $\Delta(L_0\setminus L)$ is the order complex of the poset
$(L_0\setminus L, \subseteq)$.
\end{prop}

\medskip\noindent
{\bf Proof:} As in the proof of Proposition~\ref{Sarkaria inequality} the relation (\ref{eqn:Sark-inequality}) follows from the existence of an equivariant map
\begin{equation}\label{eqn:more}
   \Delta(L_0) \stackrel{G}{\longrightarrow} \Delta(L) \ast
   \Delta(L_0\setminus L).
\end{equation}

\bigskip
The following extension of Proposition~\ref{prop:Sarkaria} is
needed in the proof of Theorem~\ref{thm:oblak}. For an orientation and an elementary introduction into the theory of diagrams of spaces the reader is referred to    \cite{Ziv98}.

\begin{prop}\label{prop:Sarkaria-Zivaljevic}{\rm (Index inequality for diagrams of spaces)}
Let $G = (\mathbb{Z}_p)^k$ be an elementary abelian $p$-group.
Suppose that $P$ is a finite (not necessarily free) $G$-poset and let $P\subseteq P_0$ be its initial, $G$-invariant
subposet. Let $P_1 = P_0\setminus P$ be the complementary subposet
of $P_0$. Assume that $\mathcal{D}_0 : P_0 \rightarrow Top$ is a
$G$-diagram of spaces with $G$-action on $\mathcal{D}$ compatible
with the action on $P_0$ and let $\mathcal{D}$ and $\mathcal{D}_1$
be the restrictions of this diagram on $P$ and $P_1$
respectively. Then,
\begin{equation}\label{eqn:S-Ziv-inequality}
{\rm Ind}_G(\| \mathcal{D}\|) \geq {\rm Ind}_G(\|
\mathcal{D}_0\|) - {\rm Ind}_G(\| \mathcal{D}_1\|) -1,
\end{equation}
where $\| \mathcal{E}\| = \mathbf{hocolim}(\mathcal{E})$ is the
homotopy colimit of the diagram $\mathcal{E}$.
\end{prop}

\medskip\noindent
{\bf Proof:}  Following a similar idea as in the proofs of Propositions~\ref{Sarkaria inequality} and \ref{prop:Sarkaria} we observe that there exists a $G$-equivariant map
\begin{equation}\label{eqn:more-2}
   \| \mathcal{D}_0\| \stackrel{G}{\longrightarrow} \| \mathcal{D}\| \ast
   \| \mathcal{D}_1\| \,.
\end{equation}
Note that in light of \cite[Proposition 3.7.]{Ziv98}  the equation (\ref{eqn:more-2}) is nothing but the equation (\ref{eqn:more}) applied to the topological poset associated to the diagram  $\mathcal{D}_0$. \hfill $\square$

\subsection{Deleted joins of $r$-unavoidable complexes}
\label{sec:del-join}

\begin{theo}\label{thm:izvor}
Suppose that $K$ is an $r$-unavoidable complex with vertices in
$[m]$. Suppose that $r=p^k$ is a prime power and let $G =
(\mathbb{Z}_p)^k$ be an elementary abelian $p$-group acting freely
on the set $[r]$. Let $K^{\ast r}_\Delta$ be the $r$-fold
($2$-wise) deleted join of $K$. Then,
\begin{equation}\label{eqn:izvor}
{\rm Ind}_G (K^{\ast r}_\Delta) \geq m - r\,.
\end{equation}
 Moreover, if $K$ is $(r,s)$-unavoidable and $G = \mathbb{Z}_p$ ($r=p$) then,
\begin{equation}\label{eqn:izvor-r-and-s}
{\rm Ind}_G (K^{\ast r}_\Delta) \geq m - r+s -1\,.
\end{equation}
  \end{theo}
\noindent
{\bf Proof of Theorem~\ref{thm:izvor}:} We apply the `Sarkaria's
inequality' (\ref{eqn:Sark-inequality})
(Proposition~\ref{prop:Sarkaria}). Let $\Delta([m])\cong
\Delta^{m-1}$ be the $(m-1)$-dimensional simplex spanned by $[m]$
and let $L_0 = \Delta([m])^{\ast r}_\Delta$. Since,
\[
L_0 = \Delta([m])^{\ast r}_\Delta \cong [r]^{\ast m},
\]
is an $(m-1)$-dimensional, $(m-2)$-connected, free $G$-complex, we
observe that ${\rm Ind}_G(L_0) = m-1$.

If $L = K^{\ast r}_\Delta$ then $L_0\setminus L$ can be described
as the set of all simplices $(A_1,\ldots , A_r)\in L_0$ such that
not all of the sets $A_i$ belong to $K$. In other words $\tau = (A_1,\ldots ,
A_r)\in L_0\setminus L$ if and only if,
\begin{equation}\label{eqn:map-phi}
\phi(\tau) = \phi(A_1,\ldots , A_r) \stackrel{def}{=} \{i\in [r]
\mid A_i\notin K\}\neq\emptyset .
\end{equation}
Observe that $\phi(\tau)\neq [r]$ since by assumption $K$ is
$r$-unavoidable. Let $\partial ([r]) = \{I\subseteq [r]\mid
\emptyset \neq I\neq [r]\}$ be the boundary poset of the simplex
$\Delta([r])$. The map,
\begin{equation}
\phi : L_0\setminus L \rightarrow \partial([r])
\end{equation}
is clearly monotone and $G$-equivariant, so it induces a
$G$-equivariant map of order complexes,
\begin{equation}
\hat{\phi} : \Delta(L_0\setminus L) \rightarrow
\Delta(\partial([r])).
\end{equation}
Since ${\rm Ind}_G(\Delta(\partial([r]))) =  {\rm
Ind}_G(S^{r-2}) = r-2$, by the monotonicity of the index
function we observe that, ${\rm Ind}_G(\Delta(L_0\setminus
L))\leq r-2$ and the inequality (\ref{eqn:izvor})  as an immediate
consequence of the Sarkaria's inequality
(\ref{eqn:Sark-inequality}).

\medskip
The inequality (\ref{eqn:izvor-r-and-s}) is established by a similar argument.
In this case the map $\phi$ (described by (\ref{eqn:map-phi})) in addition to $\phi(\tau)\neq \emptyset$ has the property that $\vert \phi(\tau)\vert\leq r-s$ for each $\tau = (A_1,\ldots , A_r)\in L_0\setminus L$. (This is an immediate consequence of the fact that the complex $K$ is $(r,s)$-unavoidable.)

\medskip
Let $P^r_s := \{I\subseteq [r]\mid
\emptyset \neq I\leq r-s\}$  be a $\subseteq$-poset. As before, there arises a $G$-map
\begin{equation}
\hat{\phi} : \Delta(L_0\setminus L) \rightarrow
\Delta(P^r_s)\,.
\end{equation}
From here we conclude that
\begin{equation}\label{eqn:conclude}
  {\rm Ind}_G(\Delta(L_0\setminus L)) \leq {\rm Ind}_G(\Delta(P^r_s)) \leq {\rm dim}(\Delta(P^r_s)) = r-s-1
\end{equation}
and the inequality (\ref{eqn:izvor-r-and-s}) follows from Sarkaria's inequality.  \hfill $\square$

\begin{cor}\label{cor:pi}
Suppose that $K$ is a simplicial complex such that ${\rm
Vert}(K)\subseteq [m]$. If $r = \pi(K)=p^k$ is a prime power and
$G=(\mathbb{Z}_p)^k$ then,
\begin{equation}\label{eqn:izvor-2}
{\rm Ind}_G (K^{\ast r}_\Delta) \geq m - \pi(K).
\end{equation}
\end{cor}

\begin{example}\label{exam:vanKampen-Flores}{\rm
The van Kampen-Flores theorem \cite[Theorem~5.1.1]{M} says that
the $n$-skeleton $\Delta^n_{2n+2}$ of the $(2n+2)$-dimensional
simplex is not embeddable in $\mathbb{R}^{2n}$. Here we deduce
this result from Theorem~\ref{thm:izvor}.

The complex $\Delta^n_{2n+2}$ is self-dual (minimally
$2$-unavoidable), hence by Theorem~\ref{thm:izvor} ${\rm
Ind}_{\mathbb{Z}_2}((\Delta^n_{2n+2})^{\ast 2}_\Delta)\geq
(2n+3)-2 = 2n+1$. If $\Delta^n_{2n+2}$ were embeddable in
$\mathbb{R}^{2n}$ then there would exist a
$\mathbb{Z}_2$-equivariant map $F : (\Delta^n_{2n+2})^{\ast
2}_\Delta \rightarrow \mathbb{R}^{4n+1}\setminus
\mathbb{R}^{2n}\simeq S^{2n}$ which contradicts the fact that
${\rm Ind}_{\mathbb{Z}_2}(S^{2n}) = 2n$. }
\end{example}

\begin{rem}
{\rm
In the proof of (\ref{eqn:izvor-r-and-s}), the stronger condition $G=\mathbb{Z}_p$ (as opposed to $G=(\mathbb{Z}_p)^k$) was needed in the proof of the inequality  ${\rm Ind}_G(\Delta(P^r_s)) \leq {\rm dim}(\Delta(P^r_s))$.
}
\end{rem}

\subsection{Deleted products of $r$-unavoidable complexes}
\label{sec:del-product}

\begin{theo}\label{thm:oblak}
Suppose that $r = p^k$ is a prime power and let $K$ be an
$r$-unavoidable complex with vertices in $[m]$. Suppose that $G=
(\mathbb{Z}_p)^k$ is a $p$-torus.  Let $K^{r}_\Delta$ and $K^{\ast r}$ be respectively  the $r$-fold ($2$-wise) deleted product (deleted join) of $K$. Then,
\begin{equation}\label{eqn:oblak}
{\rm Ind}_G (K^{r}_\Delta) \geq {\rm Ind}_G (K^{\ast r}_\Delta) -r+1 \geq m - 2r+1 \,.
\end{equation}
\end{theo}

\noindent
{\bf Proof of Theorem~\ref{thm:oblak}:} Let $P_0$ be the face poset
of the simplex $\Delta([r])$ spanned by vertices $[r]$  and let $P = \{\hat{1}\}$ be its one-element subposet where $\hat{1} = [r]$  is the maximum element of $P_0$. Let $P_1 = P_0\setminus P = P_0\setminus \{\hat{1}\}$. Define the diagram  $\mathcal{D}_0 : P_0 \rightarrow Top$ of spaces by $\mathcal{D}_0(I) = K^I_\Delta$  (for each $I\in P_0$) where $K^I_\Delta \cong K^{\vert I\vert}_\Delta$ is the $I$-deleted product of $K$ (topologized as a subspace of
$K^{[r]}$). More explicitly a function $f : I\rightarrow K$ is in
$K^{I}_\Delta$ if and only if the supports (in $\Delta([r])$) of
all points $f(i)$ are pairwise disjoint. Morever (for the inclusion $e: I\subseteq J$) the associated map $\mathcal{D}_0(e) : K^{J}_\Delta \rightarrow K^{I}_\Delta$ is the natural projection.

Let $\mathcal{D}$ and $\mathcal{D}_1$ be the restriction of this diagram to posets $P$ and $P_1$ respectively.

\medskip
It is not difficult to see \cite[Section~3.4.]{Ziv98} that the
homotopy colimit of the diagram $\mathcal{D}_0$  is the deleted join of $K$,
\begin{equation}\label{eqn:hocolim-join}
\| \mathcal{D}_0\| = \mathbf{hocolim}(\mathcal{D}_0) \cong K^{\ast
r}_\Delta.
\end{equation}
Moreover $\| \mathcal{D}\| \cong K^{r}_\Delta$ is the deleted
product of $K$.

\medskip
Let $\mathcal{C} : P_1 \rightarrow Top$ be the constant diagram
where $\mathcal{C}(I)$ is a one-element set for each $I\in P_1$.
The obvious map (morphism) of diagrams $\mathcal{D}_1 \rightarrow
\mathcal{C}$ induces a $G$-equivariant map,
\begin{equation}\label{eqn:map-diags}
\| \mathcal{D}_1\| = \mathbf{hocolim}(\mathcal{D}_1) \rightarrow
\mathbf{hocolim}(\mathcal{C})\cong \Delta(P_1) \cong S^{r-2}.
\end{equation}
By the monotonicity of the index function we conclude from
(\ref{eqn:map-diags}) that ${\rm Ind}_G(\| \mathcal{D}_1\|)
\leq r-2$. By Theorem~\ref{thm:izvor} ${\rm Ind}_G (K^{\ast
r}_\Delta) \geq m - r$. Finally by the index inequality
\ref{prop:Sarkaria-Zivaljevic},
\[
{\rm Ind}_G (K^{r}_\Delta) \geq  (m-r) - (r-2) - 1 = m - 2r
+1\,.
\]
 \hfill $\square$
\begin{rem}{\rm
  The inequality (\ref{eqn:oblak}) seems to indicate that in many applications the deleted product is at least as useful as the deleted join.  }
\end{rem}

\section{ Almost $r$-embeddings and\\ $r$-unavoidable complexes }

Suppose that $K$ is a finite simplicial complex and let $X = \vert
K\vert$ be the underlying topological space. A map $f :
X\rightarrow \mathbb{R}^d$ is an embedding if it is $1$--$1$. We
say that $f : \vert K\vert\rightarrow \mathbb{R}^d$ is an {\em
almost embedding} if $f(\Delta_1)\cap f(\Delta_2)=\emptyset$ for
each pair $\Delta_1, \Delta_2$ of vertex disjoint simplices in
$K$. In other words a map is an embedding (almost embedding) if it
does not have a $2$-fold (``strong'' $2$-fold) point. It is natural to
extend this definition to the general case of maps which do not
admit (strong) $r$-fold points.

\begin{defin}\label{def:r-embeddings}
A map $X\rightarrow \mathbb{R}^d$ is an $r$-embedding if it has no
$r$-fold points, i.e.\ there are at most $(r-1)$ points in the
pre-image $f^{-1}(y)$ for each $y\in \mathbb{R}^d$. A map $f :
K\rightarrow \mathbb{R}^d$ without {\em ``strong'' $r$-fold points} is
called an {\em almost $r$-embedding}. By
definition $y\in \mathbb{R}^d$ is a {\em strong $r$-fold point} if there
exist vertex disjoint simplices $\Delta_1,\ldots , \Delta_r$ in
$K$ such that $y\in f(\Delta_i)$ for each $i=1,\ldots, r$.
\end{defin}

\medskip
Isaac Mabillard and Uli Wagner opened a new chapter of the theory
of almost $r$-embeddings by introducing and developing a version
of `Whitney trick' for eliminating strong $r$-fold points. The
following theorem was originally announced in \cite{MW14} with the
complete presentation given in \cite{MW15}, see also \cite{amsw,
MW16} for the subsequent development.

\begin{theo}\label{thm:Ma-Wa}
{\rm (I. Mabillard, U.~Wagner \cite{MW14, MW15})} Suppose that
$r\geq 2, k\geq 3$, and let $K$ be a simplicial complex of
dimension $(r-1)k$. Then the following statements are equivalent:

\begin{enumerate}
 \item[{\rm (i)}] There exists an $S_r$-equivariant map $F : K^{\times r}_\Delta\rightarrow
S(W_r^{\oplus rk})$.
 \item[{\rm (ii)}] There exists a continuous map $f : K\rightarrow \mathbb{R}^{rk}$
such that $f(\sigma_1)\cap\cdots\cap f(\sigma_r)=\emptyset$ for
each collection of pairwise disjoint faces
$\sigma_1,\ldots,\sigma_r$ of $K$.
\end{enumerate}
\end{theo}

The following result of Murad \"{O}zaydin was never formally
published but its preprint from 1987 was widely circulated.

\begin{theo}{\rm (M.~\"{O}zaydin \cite{oz87})}\label{thm:Ozaydin}
Assume $d \geq 1$ and $r \geq 2$. Let $S_r$ be the symmetric group
and let $E_{S_r}^{d(r-1)}$ be a $[d(r-1)-1]$-connected, free
$S_r$-simplicial complex. Then there exists an $S_r$-equivariant
map,
\begin{equation}\label{eqn:Ozaydin}
f : E_{S_r}^{d(r-1)} \stackrel{S_r}{\longrightarrow} S^{d(r-1)-1}
\end{equation}
if and only if $r$ is {\em not} a prime power.
\end{theo}

Theorem~\ref{thm:Ma-Wa} and Theorem~\ref{thm:Ozaydin} together
explain why it is in general necessary to assume that $r=p^k$ is a
prime power in the Van Kampen-Flores type results. This will be
one of our standard assumptions throughout all of
Section~\ref{sec:tacit}.

\subsection{The general Tverberg-Van Kampen-Flores problem  }
\label{sec:tacit}

\begin{defin}\label{def:study-embeddings}
(Cont.\ of Definition~\ref{def:r-embeddings}) We say that a
simplicial complex $K$ is {\em almost $r$-non-embeddable in
$\mathbb{R}^d$} (or (almost) $(r,d)$-non-embeddable) if there does
not exist a map $f : \vert K\vert \rightarrow \mathbb{R}^d$
without strong $r$-fold points. The general Tverberg-Van
Kampen-Flores problem is to find interesting examples of almost
$(r,d)$-non-embeddable complexes and to study (characterize) the
class $TvKF(r,d)$ of all simplicial complexes which are almost
$r$-non-embeddable in $\mathbb{R}^d$.
\end{defin}

\medskip
The following theorem of G.~Schild illustrates the importance of
$2$-unavoidable (self-dual) complexes for the instance $r=2$  of
the general Tverberg-Van Kampen-Flores problem. In the special
case when $K_i$ are Van Kampen-Flores complexes
$\Delta^{n_i}_{2n_i+2}$ (Example~\ref{exam:vanKampen-Flores}) this
result was proved by B.~Gr\"{u}nbaum, so Theorem~\ref{thm:Schild}
is sometimes referred to as the Van
Kampen-Flores-Gr\"{u}nbaum-Schild non-embedding theorem. Sergey
Melikhov \cite{mel11} discovered an interesting connection of this
result (and its relatives) with the so called {\em dichotomial
cell complexes} (dichotomial spheres).

\begin{theo}{\rm (G.~Schild \cite{sc93}, B.~Gr\"{u}nbaum \cite{gr69})}\label{thm:Schild}
Let $K = K_1\ast\ldots\ast K_s$ where each $K_i$ is a self-dual
subcomplex of the simplex $\Delta^{m_i-1} = \Delta([m_i])$ spanned
by $m_i$ vertices. Then $K$ is not embeddable in $\mathbb{R}^d$
where,
\begin{equation}\label{eqn:Schild}
d \leq m_1+\ldots +m_s-s-2.
\end{equation}
\end{theo}

It is clear that we can replace in Theorem~\ref{thm:Schild}
self-dual complexes by $2$-unavoidable complexes. Indeed, the
$2$-unavoidable complexes which are minimal (in the sense that
each proper subcomplex $K'\subseteq K$ is not $2$-unavoidable) are
precisely the self-dual complexes. For this reason the following
theorem is a direct generalization of Theorem~\ref{thm:Schild} to
the case of $r$-unavoidable complexes.

\begin{theo}\label{thm:Dule-Sinisa-Rade}
Suppose that $r=p^k$ is a prime power. Let $K = K_1\ast\ldots\ast
K_s$ where each $K_i$ is an $r$-unavoidable subcomplex of the
simplex $\Delta^{m_i-1} = \Delta([m_i])$ spanned by $m_i$
vertices. Then $K$ is almost $r$-non-embeddable in
$\mathbb{R}^d$ if the dimension $d$ satisfies the inequality,
\begin{equation}\label{eqn:Dule-Sinisa-Rade}
(r-1)(d+s+1)+1 \leq m_1+\ldots +m_s.
\end{equation}
\end{theo}

\medskip\noindent
{\bf Proof:} The operation of the $r$-th deleted join commutes
with the standard joins (see the proof of Lemma~5.5.2 in \cite{M})
so there is an isomorphism,
 \begin{equation}\label{eqn:join-del-join}
K^{\ast r}_\Delta = (K_1\ast\cdots\ast K_s)^{\ast r}_\Delta \cong
(K_1)^{\ast r}_{\Delta}\ast\cdots\ast (K_s)^{\ast r}_{\Delta}\,.
 \end{equation}
If there exists an almost $r$-embedding  $f : K \rightarrow
\mathbb{R}^d$ then there is an associated $S_r$-equivariant map,
\begin{equation}\label{eqn:eq-map-glavna-tm}
K^{\ast r}_\Delta \stackrel{F}{\longrightarrow}
(\mathbb{R}^d)^{\ast r}_\Delta \hookrightarrow
\mathbb{R}^{dr+r-1}\setminus \mathbb{R}^d \simeq S^{(r-1)(d+1)-1} \,.
\end{equation}
Since ${\rm Ind}_G(S^{(r-1)(d+1)-1}) = N := (r-1)(d+1)-1$, in light of the monotonicity property of the index (Proposition~\ref{prop:ind-sharp}, part (3)), it is sufficient
to check the inequality,

\begin{equation}\label{eqn:when-ind}
{\rm  Ind}_G(K^{\ast r}_\Delta) \geq N+1 = (r-1)(d+1)\,.
\end{equation}
We already know (Theorem~\ref{thm:izvor}) that ${\rm
Ind}_G((K_i)^{\ast r}_\Delta)\geq m_i - r$ for each
$i=1,\ldots, s$. Hence, in light of (\ref{eqn:join-del-join}), it
may be tempting to apply a result similar to the inequality (3) in
Proposition~\ref{prop:ind-sharp}. Unfortunately this inequality
points in the `wrong direction' and in general cannot be improved to
the equality. For this reason we go back to the idea of the proof
of Theorem~\ref{thm:izvor} and use directly the
Proposition~\ref{prop:Sarkaria} (Sarkaria's inequality).

\medskip
By adding a `slack index' we can see $K_j$ as a complex with
vertices in $[m_j]\times\{j\}\subseteq \mathbb{N}_0\times
\mathbb{N}_0$. Then (in light of (\ref{eqn:join-del-join})) a
simplex $A\in K^{\ast r}_\Delta$ can be described as a collection
$A = (A_{i,j})_{(i,j)\in [r]\times [s]}$ where,

\begin{enumerate}
 \item[(1)] $A_{i,j}\subseteq [m_j]\times \{j\}$ is a simplex in
 $K_j$,
 \item[(2)] $(\forall j=1,\ldots, s)(\forall i_1\neq i_2) \quad A_{i_1,j}\cap
 A_{i_2,j}=\emptyset$,
 \item[(3)] $\bigcup_{(i,j)\in [r]\times [s]} A_{i,j}\neq\emptyset $.
\end{enumerate}
In other words $K$ is naturally interpreted as a subcomplex of the
simplex $\Delta(S)\cong \Delta([m])\cong \Delta^{m-1}$ spanned by
$S$ where $S=\cup_{j=1}^s~[m_j]\times\{j\}$ and $m = m_1+\ldots +
m_s$.

\medskip
Let $L_0 = \Delta(S)^{\ast r}_\Delta \cong [r]^{\ast m}$ which
implies the inequality ${\rm  Ind}_G(L_0)\geq m-1$. Let $L
= K^{\ast r}_\Delta$. Then $L_0\setminus L$ is the collection $A =
(A_1,\ldots, A_r)\in L_0 = \Delta(S)^{\ast r}_\Delta$ such that
$A_i\notin K$ for at least one $i\in [r]$. Let $A_{i,j} = A_i\cap
([m_j]\times\{j\})$. Then, $A = (A_{i,j})_{(i,j)\in [r]\times
[s]}$ is in $L_0\setminus L$ if and only if,

\begin{enumerate}
 \item[(1)] $(\exists i\in [r])(\exists j\in [s]) \quad
 A_{i,j}\notin K_j$,
 \item[(2)] $(\forall j\in [s])(\forall i_1\neq i_2) \quad A_{i_1,j}\cap
 A_{i_2,j}=\emptyset$,
 \item[(3)] $\bigcup_{(i,j)\in [r]\times [s]} A_{i,j}\neq\emptyset $.
\end{enumerate}

\medskip
For a given $A\in K^{\ast r}_\Delta$, let $\phi_j(A) = \{i\in
[r]\mid A_{i,j}\notin K_j\}$. Since $K_j$ is $r$-unavoidable,
$\phi_j(A)\in \partial\Delta([r])\cong S^{r-2}$. Moreover, there
is monotone map (of $\subseteq$-posets),
\begin{equation}
\phi = \phi_1\uplus\ldots \uplus \phi_s: L_0\setminus L
\longrightarrow (\partial\Delta([r]))^{\ast s}
\end{equation}
which, in light of the isomorphism $(\partial\Delta([r]))^{\ast
s}\cong (S^{r-2})^{\ast s}$, leads to the inequality ${\rm
Ind}_G(\Delta(L_0\setminus L))\leq (r-1)s-1$. Together with
the Sarkaria's inequality this implies,
\begin{equation}\label{eqn:za-kraj}
{\rm  Ind}_G(L) \geq m-1-[(r-1)s-1]-1 = m-(r-1)s-1\,.
\end{equation}
This completes the proof of the theorem since the inequality
(\ref{eqn:when-ind}) is an immediate consequence of
(\ref{eqn:Dule-Sinisa-Rade}) and (\ref{eqn:za-kraj}). \hfill
$\square$

\begin{example}\label{exam:type-B}{\rm (\cite{VZ94})
The $0$-dimensional complex $L = [5]\subseteq 2^{[5]}$ is clearly
$3$-unavoidable. It follows from
Theorem~\ref{thm:Dule-Sinisa-Rade} that the complex $K = [5]\ast
[5]\ast [5]$ is almost $3$-non- embeddable in $\mathbb{R}^3$.
Informally speaking we claim that for each constellation of red,
blue and white stars in the outer space there exist three
intersecting, vertex disjoint triangles with vertices formed by
the stars of different color. This result is an instance of the
{\em type~B} colored Tverberg theorem of Vre\' cica and \v
Zivaljevi\' c, \cite{VZ94, Z04}. As shown in \cite{VZ94} the
complex $[4]\ast [n] \ast [n]$ is almost $3$-embeddable into
$\mathbb{R}^3$ for all $n\in \mathbb{N}$. }
\end{example}

\section{Applications}

In this section we recast Theorem~\ref{thm:Dule-Sinisa-Rade} in the language of cooperative game theory. There are several reasons why this reformulation may be interesting. First of all we show how, by using {\em threshold complexes}, we can automatically construct  $r$-unavoidable complexes and, in light of Theorem~\ref{thm:Dule-Sinisa-Rade}, generate numerous examples of almost $r$-non-embeddable complexes.
This connection of topology and cooperative game theory, albeit at first sight somewhat superficial, may nevertheless provide an interesting interface for exchange of ideas, methods  and questions from the two apparently distant mathematical fields.

\subsection{Simple games and threshold complexes}

A {\em simple game} is a mathematical concept used in cooperative game theory  to
describe the distribution of power among coalitions of players. The simplest definition of a simple game (from the viewpoint of a combinatorial topologist) says that a family of sets $F\subseteq 2^{[n]}$ is a simple game if and only if the complement $K := 2^{[n]}\setminus F$ is a simplicial complex. However, simple games have appeared in a variety of other mathematical
contexts under various names, including the following:

\begin{enumerate}
  \item threshold functions,
  \item boolean or switching functions,
  \item hypergraphs,
  \item coherent structures,
  \item Sperner systems,
  \item abstract simplicial complexes.
\end{enumerate}

\begin{defin}
  A simple game $F\subset 2^{[n]}$ is called a {\em weighted majority game} if there exists a non-negative weight distribution $w = (w_1,\dots, w_n)\in \mathbb{R}^n_+$ and a real number $q$, called quota, such that
  \[
         X \in F  \qquad   \Leftrightarrow  \qquad    w(X):= \sum_{i\in X} w_i  > q  \,.
  \]
  Such a simple game is often recorded as the pair $[q; w] = [q; w_1,\dots, w_n]$.
  The complementary simplicial complex $K_{w\leq q} := 2^{[n]}\setminus F$ is referred to as the {\em threshold complex} with weight distribution $w$ and quota $q$.
\end{defin}

\medskip
The following elementary lemma explains why we are interested in threshold complexes.

\begin{lema}\label{lem:threshold}
 Suppose that the total weight of $w$ is  ${\tilde{w}} := w([n]) = w_1+\dots + w_n $.
Then the threshold complex,
\[
    K_{w \leq {\tilde{w}}/{r}} := \{X\subset [n]  \mid w(X)\leq {\tilde{w}}/{r}\}
\]
is $r$-unavoidable.
 \end{lema}

\subsection{Theorem~\ref{thm:Dule-Sinisa-Rade} revisited}

The following reformulation of Theorem~\ref{thm:Dule-Sinisa-Rade} provides an efficient algorithm for generating almost $r$-non-embeddable complexes.  The reader is invited to test and apply  this algorithm to non-planarity of graphs (Kuratowski graphs), Van Kampen-Flores theorem, and other results of Tverberg-Van Kampen-Flores type.

\begin{enumerate}
\item Choose $r$, quota $q = \tilde{w}/r$ and the number $m = m_1$ of
players.

\item Choose non-negative weights $w_1,\dots, w_{m_1}$ and construct the associated simple game $[q; w] = [q; w_1,\dots, w_{m_1}]$ where $q:=\tilde{w}/r$.

\item Determine the associated $r$-unavoidable threshold complex
$K_1 := K_{w\leq \tilde{w}/r}$.

\item Repeat the procedure $s\geq 1$ times, possibly changing the
number of players (and the corresponding weights), using the quota $q=\tilde{w}/r $.

\item Record the associated threshold complexes $K_i\subset 2^{[m_i]}$.

\item  Let $K = K_1\ast\dots\ast K_s$ be the associated join.

\item Find $d$ from the equation,
\[
   (r-1)(d+s+1)+1 = m_1+\dots+ m_s.
\]
\end{enumerate}

Then, $K$ is  almost $r$-non-embeddable in
$\mathbb{R}^d$.

\begin{rem}{\rm
  There exist classes of $r$-unavoidable complexes which are not threshold complexes and which are unavoidable for deeper reasons. Examples include joins of minimal triangulations of manifolds which `look like a projective plane', see
  \cite[Section 7]{bestiary}. These complexes can be also used as an input for Theorem~\ref{thm:oblak} and the algorithm for generating almost $r$-non-embeddable complexes, described in this section.  }
\end{rem}


\begin{thebibliography}{10000}


\bibitem[AMSW]{amsw} S. Avvakumov, I. Mabillard, A. Skopenkov, U.
Wagner. Eliminating higher-multiplicity intersections, III.
Codimension 2, arXiv:1511.03501 [math.GT].

 \bibitem[Bar]{Bar}
        T. Bartsch. \textit{Topological Methods for Variational Problems with Symmetries}.
        Lecture Notes in Mathematics 1560, Springer-Verlag Berlin Heidelberg (1993).
\bibitem[BMZ]{bmz}
P.V.M.~Blagojevi{\'c}, B.~Matschke, G.M.~Ziegler.
 Optimal bounds for the colored Tverberg problem. {\em Advances in Math.}, 226 (2011), 5198--5215.  arXiv:0910.4987v2 [math.CO].

\bibitem[BFZ]{bfz}
P.V.M.~Blagojevi{\'c}, F.~Frick, G.M.~Ziegler.
 Tverberg plus constraints.
 {\em B. London Math. Soc.}, {\textbf{46}}, (2014), 953--967.

\bibitem[BFZ-2]{bfz2}
P.V.M.~Blagojevi{\'c}, F.~Frick, G.M.~Ziegler.
 Barycenters of polytope skeleta and counterexamples to the topological tverberg conjecture, via constraints. {\em J. Europ. Math. Soc. (JEMS)} (to appear), 	arXiv:1510.07984 [math.CO].

\bibitem[ClPu]{CP1}
        M.~Clapp, D.~Puppe. \textit{Critical point theory with symmetries}.
        J. reine angew. Math., {\textbf{418}}, (1991), 1-29.

\bibitem[Die87]{dieck}
T.~tom Dieck. {\em Transformation Groups}. Vol.~8 of \emph{Studies
in Mathematics}. Walter de Gruyter, Berlin, 1987.

\bibitem[F]{Frick} F.~Frick.
 Counterexamples to the topological Tverberg conjecture. {\em Oberwolfach Reports} 12(1) (2015), 318--321.  arXiv:1502.00947v1 [math.CO].

\bibitem[Ga-Pa]{ga-pa} P.~Galashin, G.~Panina. Simple game induced
manifolds.  {\em J. Knot Theory Ramif.} 25(12), (2016).   arXiv:1311.6966 [math.GT].

\bibitem[G10]{g10} M.~Gromov. Singularities, expanders and topology of
maps, Part 2: From combinatorics to topology via algebraic
isoperimetry. \textit{Geom. Funct. Anal.}, {\textbf{20}}, (2010),
416–-526.

\bibitem[Gr69]{gr69} B.~Gr\"{u}nbaum. Imbeddings of simplicial
complexes. \textit{Comment. Math. Helv.}, {\textbf{44}}, (1969),
502--513.

\bibitem[J-M]{j-m} J.~Jezierski, W.Marzantowicz. \textit{Homotopy Methods in Topological Fixed and Periodic Points Theory}. Volume 3 of Topological Fixed Point Theory and Its Applications, Springer  2006.

\bibitem[JVZ-1]{jvz-1}
D.~Joji{\'c}, S.T.~Vre{\'c}ica, R.T.~{\v Z}ivaljevi{\'c}.
 Multiple chessboard complexes and the colored Tverberg
problem. {\em J. Combin. Theory Ser. A} 
 145,  2017, 400--425.  arXiv:1412.0386 [math.CO].

\bibitem[JVZ-2]{jvz-2}
D.~Joji{\'c}, S.T.~Vre{\'c}ica, R.T.~{\v Z}ivaljevi{\'c}.
 Symmetric multiple chessboard complexes and a new theorem of Tverberg
type.   {\em J. Algebraic Combin.}, (2017), 46,  15--31. 
arXiv:1502.05290v2 [math.CO].

\bibitem[JVZ-3]{jvz-3}
D.~Joji{\'c}, S.T.~Vre{\'c}ica, R.T.~{\v Z}ivaljevi{\'c}. Topology
and combinatorics of `unavoidable complexes'. arXiv:1603.08472
[math.AT].

\bibitem[JJTVZ]{bestiary}
M.~Jeli\'{c}, D.~Joji{\'c}, M.~Timotijevi\'{c}, S.T.~Vre{\'c}ica,
R.T.~{\v Z}ivaljevi{\'c}.  Combinatorics of unavoidable complexes.
arXiv:1612.09487 [math.CO].

\bibitem[MW14]{MW14} I.~Mabillard and U.~Wagner. Eliminating Tverberg
points, I; An analogue of the Whitney trick. In {\em Proc. 30th Ann.
Symp. on Computational Geometry}, (2014),  171--180.


\bibitem[MW15]{MW15} I.~Mabillard, U.~Wagner.  Eliminating Higher-Multiplicity
Intersections, I. A Whitney Trick for Tverberg-Type Problems.
Preprint, 46 pages, arXiv:1508.02349.

\bibitem[MW16]{MW16} I.~Mabillard, U.~Wagner. Eliminating Higher-Multiplicity Intersections, II. The Deleted Product Criterion in the r-Metastable Range.
arXiv:1601.00876 [math.GT].


\bibitem[Mar]{Marzantowicz1}
W.~Marzantowicz. {\it A G-Lusternik--Schnirelman category of space
with an action of a compact Lie group}. Topology {\textbf{28}},
(1989), 403-412.

\bibitem[M03]{M} J.~Matou\v sek. \textit{Using the Borsuk-Ulam Theorem.
Lectures on Topological Methods in Combinatorics and Geometry}.
Universitext, Springer-Verlag, Heidelberg, 2003 (Corrected 2nd
printing 2008).

\bibitem[Mel11]{mel11} S.A.~Melikhov. Combinatorics of embeddings, arXiv:1103.5457v2
[math.GT].


\bibitem[vNM44]{vNM44} J.~Von Neumann. O.~Morgenstern. \textit{Theory of Games and Economic
Behavior}, Princeton University Press, 1944.

\bibitem[\"{O}z87]{oz87} M.~\"{O}zaydin. Equivariant maps for the symmetric group. Unpublished  preprint.  \url{http://minds.wisconsin.edu/handle/1793/63829}, 1987. 

\bibitem[PS07]{ps07} B.~Peleg, P.~Sudh\"{o}lter. \textit{Introduction to the Theory of
Cooperative Games}. Springer Science \& Business Media, 2007.

\bibitem[R90]{r90} K.G.~Ramamurthy. \textit{Coherent Structures and Simple
Games}, Springer Netherlands, 1990.


\bibitem[Sc93]{sc93} G. Schild. Some minimal nonembeddable complexes, \textit{Topology
Appl.} {\textbf{3}}, (1993), no. 2, 177–-185.


\bibitem[Vol96]{Vol} A. Yu. Volovikov. On a topological generalization of the
Tverberg theorem. {\em Math. Notes}, {\textbf{59}}, (3), (1996)
324--326. Translation from {\em Mat. Zametki} {\textbf{59}}, No.3,
(1996), 454-456 .

\bibitem[V\v Z94]{VZ94}  S.T.~Vre\' cica and R.~\v Zivaljevi\' c. New cases of
the colored Tverberg theorem. In H.~Barcelo and G.~Kalai, editors,
\textit{Jerusalem Combinatorics '93}, Contemporary Mathematics Vol.\
178, pp. 325--334, A.M.S.\ 1994.

\bibitem[V{\v{Z}}11]{vz11}
 S.T.~Vre{\'c}ica, R.T.~{\v{Z}}ivaljevi{\'c}.
 Chessboard complexes indomitable.
 {\em J. Combin. Theory Ser. A}, {\textbf{118}}, (7), (2011),
2157--2166.

\bibitem[Zi11]{Ziegler}  G.M.~Ziegler.  \newblock $3N$ colored points in a plane.
 {\em Notices of the A.M.S.} Vol. {\textbf{58}}, Number 4,(2011),
550--557. .

\bibitem[{\v{Z}}V92]{zv92}
R.T.~{\v{Z}}ivaljevi{\'c} and  S.T.~Vre{\'c}ica.
 The colored {T}verberg's problem and complexes of injective
  functions.
 {\em J. Combin. Theory Ser. A}, {\textbf{61}}, (2), (1992),
309--318.

\bibitem[\v Ziv98]{Ziv-96-98}
R.T.~\v Zivaljevi\' c.  User's guide to equivariant methods in
combinatorics, I and II.
\newblock {\it Publ. Inst. Math. (Beograd) (N.S.)}, (I) {\textbf{59}}, (73), (1996), 114--130,  and (II)
{\textbf{64}}, (78), (1998), 107--132.

\bibitem[\v Z98]{Ziv98} R.T.~\v Zivaljevi\'{c}. Combinatorics of
topological posets. \textit{Advances in Applied
Math.},{\textbf{21}}, (1998), 547--574.

\bibitem[\v Z04]{Z04}
R.T.~\v Zivaljevi\'{c}. Topological methods. Chapter 14 in
\textit{Handbook of Discrete and Computational Geometry}, J.E.\
Goodman, J.\ O'Rourke, eds, Chapman \& Hall/CRC 2004, 305--330.

\bibitem[\v Z17]{Z17}
R.T.~\v Zivaljevi\'{c}. Topological methods in discrete geometry.
Chapter 21 in \textit{Handbook of Discrete and Computational
Geometry}, third ed., J.E. Goodman, J. O'Rourke, and C.D. T\'{o}th,
CRC Press LLC, Boca Raton, FL, 2017.




\end{thebibliography}
\end{document}